\documentstyle{amsppt}
\tolerance 3000
\pagewidth{5.5in}
\vsize7.0in
\magnification=\magstep1
\widestnumber \key{AAAAA}
\topmatter
\author Alex Iosevich, Nets Katz, and Steen Pedersen \endauthor
\address Department of Mathematics, Georgetown University, Washington, 
DC 20057, USA \endaddress
\thanks Research supported in part by NSF grants DMS97-06825 and DMS-9801410
\endthanks
\email iosevich\@math.georgetown.edu \endemail
\address Department of Mathematics, University of Illinois at Chicago, 
Chicago, Illinois, 60607 
\endaddress
\email nets\@math.uic.edu 
\endemail
\address Department of Mathematics, Wright State University, Dayton, 
OH 45435, USA \endaddress
\email steen\@math.wright.edu
\endemail
\title Fourier bases and a distance problem of Erd\H os
\endtitle
\abstract We prove that no ball admits a non-harmonic orthogonal 
basis of exponentials. We use a combinatorial result, originally 
studied by Erd\H os, which says that the number of distances determined 
by $n$ points in ${\Bbb R}^d$ is at least 
$C_d n^{\frac{1}{d}+\epsilon_d}$, 
$\epsilon_d>0$. 
\endabstract
\endtopmatter

\head Introduction and statement of results \endhead 

\subhead Fourier bases \endsubhead Let $D$ be a domain in ${\Bbb 
R}^d$, i.e., $D$ is a Lebesgue measurable subset of $\Bbb{R}^d$ with 
finite non-zero
Lebesgue measure. We say that $D$ is a {\it spectral set} if $L^2(D)$ has 
orthogonal basis of the form 
$E_{\Lambda}
={\{e^{2 \pi i x \cdot \lambda} \}}_{\lambda \in \Lambda}$, 
where $\Lambda$ is an infinite subset of 
${\Bbb R}^d$. We shall refer to $\Lambda$ as a {\it spectrum} for $D$. 

We say that a family $D+t$, $t\in T$, of translates of a domain $D$ 
{\it tiles} $\Bbb{R}^d$ if 
$\cup_{t\in T}\left(D+t\right)$ is a partition of
$\Bbb{R}^d$ up to sets of Lebesgue measure zero.

\example{Conjecture}It has been conjectured (see 
\cite{Fug}) that a domain $D$ is a spectral set if and only if it is 
possible to tile ${\Bbb R}^d$ by a family of 
translates of $D$. 
\endexample
 
This conjecture is nowhere near resolution, even in 
dimension one. It has been the subject of recent research,
see for example \cite{JoPe2}, \cite{LaWa}, and \cite{Ped}.

In this paper we address the following special case of 
the conjecture. Let $B_d=\{x \in {\Bbb R}^d: |x| \leq 1 \}$ 
denote the unit 
ball. We prove that 
\proclaim{Theorem 1} An affine image of $D=B_{d}$, $d \ge 2$, is not a 
spectral set.
\endproclaim 
If $A$ is a (possibly unbounded) 
self-adjoint operator acting on some Hilbert space, then we
may define $\exp\left(-\sqrt{-1}A\right)$ using the Spectral Theorem. We
say that two (unbounded) self-adjoint operators $A$ and $B$ acting on 
the same
Hilbert space {\it commute} if the bounded unitary operators 
$\exp\left(-\sqrt{-1}sA\right)$ and $\exp\left(-\sqrt{-1}tB\right)$
commute for all real numbers $s$ and $t$. See, for example, \cite{ReSi}
for more details on the needed operator theory.
As an immediate consequence of \cite{Fug} and Theorem 1 we have:
\proclaim{Corollary}
There do not exist commuting self-adjoint operators $H_j$ acting on
$L^2(B_d)$ such that 
$H_jf=-\sqrt{-1}\,\partial f/\partial x_j$ for $f$ in the domain
of the unbounded operator $H_j$ and $1\leq j\leq d$. 
The derivatives $\partial /\partial x_j$ act on $L^2(B_d)$ in the 
distribution sense.
\endproclaim

In other words, there do not exist commuting self-adjoint restrictions of
the partial derivative operators 
$-\sqrt{-1}\,\partial /\partial x_j$, $j=1,\ldots,d$,
acting on $L^2(B_d)$ in the distribution sense.

The two-dimensional case of Theorem 1 was proved by Fuglede in 
\cite{Fug}. Our proof uses the following combinatorial result. 
See for example \cite{AgPa}, Theorem 12.13. 

\proclaim{Theorem 2} Let $g_d(n)$, $d \ge 2$, denote the minimum number of 
distances determined by $n$ points in ${\Bbb R}^d$. Then 
$$ g_d(n) \ge C_d n^{\frac{3}{3d-2}}. \tag*$$ 
\endproclaim 

\remark{Remark} The study of the problem addressed in Theorem 2 was 
initiated by Erd\H os. He proved that $g_2(n) \ge Cn^{\frac{1}{2}}$. 
See \cite{Erd}. Moser proved in \cite{Mos} that $g_2(n) \ge 
Cn^{\frac{2}{3}}$. More recently, Chung, Szeremedi, and Trotter proved 
that $g_2(n) \ge C\frac{n^{\frac{4}{5}}}{\log^c(n)}$ for some $c>0$. 
See \cite{CST}. Theorem 2 above is proved by induction using the $g_2(n) 
\ge Cn^{\frac{3}{4}}$ result proved by Clarkson et al. in \cite{C}. 

As the reader shall see, Theorem 1 does not require the full strength of 
Theorem 2. We just need the fact $g_d(n) \ge C_dn^{\frac{1}{d}+\epsilon}$, 
for some $\epsilon>0$.    
\endremark 

It is interesting to contrast the case of the ball with the case of 
the cube ${[0,1]}^d$. It was proved in \cite{IoPe1}, (and, 
independently, in \cite{LRW}; for $d\leq3$ this was established in 
\cite{JoPe2}), that $\Lambda$ is a spectrum for  
${[0,1]}^d$, in the sense defined above, if and only if $\Lambda$ is 
a tiling set for ${[0,1]}^d$, in the sense that 
${[0,1]}^d+\Lambda={\Bbb R}^d$ without overlaps. It follows that 
${[0,1]}^d$ has lots of spectra. The standard 
integer lattice $\Lambda={\Bbb Z}^d$ is an example, 
though there are many 
non-trivial examples as well. See \cite{IoPe1} and \cite{LaSh}. 

Our method of proof is as follows. We shall argue 
that if $B_{d}$ were a spectral set, then any 
corresponding spectrum $\Lambda$ would have the property 
$\# \{ \Lambda \cap B_d(R) \} \approx R^d$, 
where $B_d(R)$ denotes a ball of 
radius $R$ and $f(R)\approx g(R)$ means that there exist constants
$c\leq C$ so that $c\,f(R)\leq g(R)\leq C\, f(R)$ for $R$ sufficiently
large. 
On the other hand, we will show that the number of distinct 
distances between the elements of $\{ \Lambda \cap B_d(R) \}$ is 
$\approx R$. Theorem 2 implies that if $R$ is sufficiently large, 
this is not possible. 

Kolountzakis (\cite{Kol}) recently proved that if $D$ is any convex 
non-symmetric domain in ${\Bbb R}^d$, then $D$ is not a spectral set. 
Theorem 1 is a step in the direction of proving 
that if $D$ is a convex domain such that 
$\partial D$ has at least one point where the Gaussian curvature 
does not vanish, then $D$ is not a spectral set. 
This, in its turn, would be a step towards proving 
the conjecture of Fuglede mentioned above. 
\vskip.125in 

\head Orthogonality \endhead 

For a domain $D$ let 
$$ Z_{D}=\{\xi\in\Bbb{R}:\widehat{\chi}_{D}(\xi)=0\}.$$
Consider a set of exponentials $E_{\Lambda}$. 
Observe that 
$$\widehat{\chi}_{D}(\lambda-\lambda')
   =\int_{D}e_{\lambda}(x)\overline{e_{\lambda'}(x)}\,dx.$$
It follows that the exponentials 
$E_{\Lambda}$ are orthogonal in $L^{2}(D)$ iff 
$$\Lambda-\Lambda\subseteqq Z_{D}\cup\{0\}.$$
\proclaim{Proposition 1}
If $E_{\Lambda}$ is an orthogonal subset of $L^{2}(D)$ then 
there exists a
constant $C$ depending only on $D$ such that 
$$\#\left(\Lambda\cap B_{d}(R)\right)\leq C\,R^{d}$$
for any ball $B_{d}(R)$ of radius $R$ in $\Bbb{R}^{d}$.
\endproclaim
\demo{Proof}
   Since $\widehat{\chi}_{D}$ is continuous and 
   $\widehat{\chi}_{D}(0)=|D|$
   it follows that 
   $$\inf\{|\xi|:\widehat{\chi}_{D}(\xi)=0\}=r>0.$$ 
   If $\xi_{1}$,
   $\ldots$, $\xi_{n}$ are in $\Lambda\cap B_{d}(R)$ then the balls 
   $B(\xi_{j},r/2)$ are disjoint and contained in $B_{d}(R+r/2)$. 
   Since $r$
   only depends on $D$ the desired inequality follows.
\enddemo

To study the exact possibilities for sets $\Lambda$ so that 
$E_{\Lambda}$ is
orthogonal it is of interest to us to compute the set $Z_{D}$. 
We will without loss of generality assume that $0\in\Lambda$. 
We again compare the sets $Z_{D}$ for the cases where $D$ is the 
cube and the ball. 

Let $Q_{d}=[0,1]^d$ be the cube in $\Bbb{R}^d$.
The zero set $Z_Q$ for $\widehat{\chi}_Q$ is the union of the
hyperplanes $\{x\in\Bbb{R}^d:x_i=z\}$, where 
the union is taken over $1\leq i \leq d$, and over all 
non-zero integers $z$.

Let $B_{d}=\{x\in\Bbb{R}^d: \|x\|\leq 1\}$ be the unit ball in
$\Bbb{R}^d$. The zero set $Z_{B_{d}}$ for $\widehat{\chi}_B$ is the 
union of the spheres
$\{x\in\Bbb{R}^d:\|x\|=r\}$, where the union is over all the positive 
roots $r$ of an appropriate Bessel
function. 

For the cube $Q_{d}$ it is easy to find a large set 
$\Lambda\subseteqq Z_{Q_{d}}\cup\{0\}$ so that 
$\Lambda-\Lambda \subseteqq Z_{Q_{d}}\cup\{0\}$. 
For example, we may take 
$\Lambda=\Bbb{Z}^{d}$. In the case of the ball $B_{d}$,
 we will show that only 
relatively small sets $\Lambda\subseteqq Z_{B_{d}}\cup\{0\}$ satisfy
$\Lambda-\Lambda \subseteqq Z_{B_{d}}\cup\{0\}$.
\vskip.125in 

\head Proof of Theorem 1 \endhead 

We shall need the following result.
\proclaim{Theorem 3} Suppose that $D$ is a spectral set and that 
$\Lambda$ is a spectrum for $D$ 
in the sense defined above, where $D$ is a bounded domain. There 
exists an $r>0$ so that any ball of radius $r$ contains at least one 
point from $\Lambda$. 
\endproclaim 
\demo{Proof}
This is a special case of \cite{IoPe2}. See also \cite{Beu}, 
\cite{Lan}, and \cite{GrRa}.
\enddemo
It is a consequence of Theorem 3 that if $D$ is a spectral set then 
there exists a constant $C>0$ such that if $\Lambda$ is a spectrum for 
$D$ then 
$\# \{ \Lambda \cap B_d(R) \}\ge C\, R^{d}$
for any ball $B_d(R)$ of radius $R$ provided that
$R$ is sufficiently large. Combining this with Proposition 1 we see 
that $\# \{ \Lambda \cap B_d(R) \}\approx R^{d}$.

Suppose $\Lambda$ is a spectrum for the unit ball $B_d$ centered 
at the origin in $\Bbb{R}^d$. 
Let $B_d(R)$ be a ball of radius $R$. 
Since $\# \{ \Lambda \cap B_d(R) \}\approx R^{d}$
it follows from Theorem 2 that 
$$\#\{|\lambda-\lambda'|:
      \lambda,\lambda'\in\Lambda\cap B_d(R)\}
   \geq C\, R^{\frac{3d}{3d-2}}.\tag**$$
Now, since $\widehat{\chi}_{B_d}$ is an analytic radial function, it follows
that if $f$ is given by $f(|\xi|)=\widehat{\chi}_{B_d}(\xi)$, 
then the number of
zeros of $f$ in the interval $[-R,R]$ is bounded above by a multiple of
$R$. In fact an explicit calculation shows that 
$\widehat{\chi}_{B_d}(\xi)
   =|\xi|^{\frac{d}{2}}J_{\frac{d}{2}}(2\pi|\xi|)$, where $J_{\nu}$
denotes the usual Bessel function of order $\nu$. See, for example,
\cite{BCT, p. 265}.

If $\lambda$, $\lambda'\in\Lambda$ then 
$$f(|\lambda-\lambda'|)= \widehat{\chi}_{B_d}(\lambda-\lambda')=0.$$
Combining the upper bound on the number of zeros of $f$ in $[-R,R]$
with the lower bound (**) we derived from Theorem 2 above we have
$$C'\, R\geq \#\{|\lambda-\lambda'|:
      \lambda,\lambda'\in\Lambda\cap B_d(R)\}
   \geq C\, R^{\frac{3d}{3d-2}}.$$
Since $1<\frac{3d}{3d-2}$ this leads to a contradiction by choosing $R$
sufficiently large.
This completes the proof of Theorem 1. 

\vskip.25in  
\head References \endhead

\ref \key AgPa \by P. Agarwal and J. Pach \book Combinatorial Geometry
\publ Wiley-Interscience Series \yr 1995 \endref

\ref \key Beu \by A. Beurling \paper Local harmonic analysis with some 
applications to differential operators \yr 1966 \jour Some Recent 
Advances 
in the Basic Sciences, Academic Press \vol 1 \endref 

\ref \key BCT \by L. Brandolini, L. Colzani, and G. Travaglini \paper 
Average decay of Fourier transforms and integer points in polyhedra 
\jour Ark. Mat. \vol 35 \yr 1997 \pages 253-275 \endref

\ref \key CST \by F. Chung, E. Szeremedi, and W. Trotter \paper 
The number of distinct distances determined by a set of points in the
Euclidean plane \jour Discrete and Computational Geometry \yr 1992 \vol 7
\endref 

\ref \key C \by K. Clarkson, H. Edelsbrunner, L. Guibas, M. Sharir, and
E. Welzl \paper Combinatorial complexity bounds for for arrangements of 
curves and surfaces \jour Discrete and Computational Geometry \vol 5 
\yr 1990 \endref 

\ref \key Erd \by P. Erd\H os \paper On sets of distances of $n$ points 
\jour American Mathematical Monthly \vol 53 \yr 1946 \endref  

\ref \key Fug \by B. Fuglede \paper Commuting self-adjoint partial 
differential operators and a group theoretic problem \jour J. Funct. 
Anal. 
\yr 1974 \vol 16 \pages 101-121 \endref 

\ref \key GrRa \by K. Gr\"{o}chenig and H. Razafinjatovo \paper On 
Landau's 
necessary density conditions for sampling and interpolation of 
band-limited
functions \jour  J. London. Math. Soc. \vol 54 \pages 557-565 \yr 
1996 
\endref 

\ref \key IoPe1 \by A. Iosevich and S. Pedersen \paper Spectral and 
tiling 
properties of the unit cube \jour Internat. Math Reseach Notices 
\vol 16 \pages 819-828\yr 1998 
\endref 

\ref \key IoPe2 \by A. Iosevich and S. Pedersen \paper How large are 
the spectral gaps? 
\jour Pacific J. Math. (to appear) \vol \pages \yr 1998
\endref

\ref \key JoPe1 \by P. E. T. Jorgensen and S. Pedersen \paper Spectral 
pairs in Cartesian coordinates
\jour J. Fourier Anal. Appl. (to appear) \vol \pages \yr 1998
\endref

\ref \key JoPe2 \by P. E. T. Jorgensen and S. Pedersen \paper Orthogonal 
harmonic analysis of fractal measures
\jour ERA Amer. Math. Soc. \vol 4 \pages 35-42 \yr 1998
\endref

\ref \key Kol \by M. Kolountzakis \paper Non-symmetric convex domains 
have no basis 
of exponentials \jour (preprint) \vol \pages \yr 1999 \endref 

\ref \key LRW \by J. Lagarias, J. Reed, and Y. Wang \paper 
Orthonormal 
bases of exponentials for the $n$-cube \jour (preprint) 
\vol \pages \yr 1998 
\endref

\ref \key LaSh \by J. Lagarias and P. Shor \paper 
Keller's cube tiling conjecture is false in high dimensions
\jour Bull. Amer. Math. Soc.
\vol 27 \pages 279-283 \yr 1992 
\endref 

\ref \key LaWa \by J. Lagarias and Y. Wang \paper 
Spectral sets and factorizations of finite abelian groups
\jour J. Funct. Anal.
\vol 145 \pages 73-98 \yr 1997 
\endref 

\ref \key Lan \by H. Landau \paper Necessary density conditions for 
sampling and interpolation of certain entire functions \jour Acta 
Math. 
\vol 117 \pages 37-52 \yr 1967 \endref 

\ref \key Mon \by H. Montgomery \paper Ten lectures on the interface 
between
analytic number theory and harmonic analysis \jour CBMS Regional 
Conference 
Series in Mathematics \yr 1994 \endref \

\ref \key Mos \by L. Moser \paper On different distances determined by 
$n$ points \jour American Mathematical Monthly \yr 1952 \vol 59 \endref 

\ref \key Ped \by S. Pedersen \paper 
Spectral sets whose spectrum is a lattice with a base
\jour J. Funct. Anal.
\vol 141 \pages 496-509 \yr 1992 
\endref 

\ref \key ReSi \by M. Reed and B. Simon 
\book Modern Mathematical Physics, Vol. I: Functional Analysis
\publ Academic Press \yr 1972 \endref

\enddocument